\documentclass{amsart}
\usepackage{amscd,amssymb}

\theoremstyle{plain}
\newtheorem{thm}{Theorem}[section]
\newtheorem{lem}[thm]{Lemma}

\newtheorem{cor}[thm]{Corollary}
\newtheorem{prop}[thm]{Proposition}

\theoremstyle{definition}

\theoremstyle{remark}

\numberwithin{equation}{section}

\begin{document}
\title{$\aleph_0$-categorical groups and their completions}
\author{ Aleksander Ivanov}

\maketitle

\centerline{Institute of Mathematics, University of Wroc{\l}aw} 
\centerline{pl.Grunwaldzki 2/4, 50-384 Wroc{\l}aw, Poland }
\centerline{E-mail: ivanov@math.uni.wroc.pl}

\bigskip 

\bigskip

\begin{quote}
ABSTRACT
\footnote{{\bf 2000 Mathematics Subject Classsification}: 
20A05, 20F18, 03C45, 20F65}
\footnote{The research was supported by KBN grant 1 P03A 025 28}
We embed a countably categorical group $G$ into 
a locally compact group $\overline{G}$ with a non-trivial 
topology and study how topological properties of $\overline{G}$ 
are connected with the structure of definable subgroups of $G$  
. 
\end{quote}
\bigskip

\section{Introduction}
Using the paper of Belyaev \cite{belyaev} we show 
that any infinite countably categorical group $G$  
has an action on some countable set $X$ which induces 
an embedding of $G$ into $Sym(X)$ such that the closure 
of $G$ in $Sym(X)$ is a locally compact subgroup 
$\overline{G}$ with non-trivial topology.  
We call the group $\overline{G}$ the {\em completion} 
of $G$ (with respect to the action). 
We will see that when $G$ is residually finite, 
the construction of the action $(G,X)$ has the property 
that the completion $\overline{G}$ is the profinite 
completion of $G$. \parskip0pt 

This provides new tools in the subject. 
In particular we show that there are strong 
connections between the lattice of 0-definable subgroups 
of $G$ and topological properties of $\overline{G}$. 
In particular the cofinality of $\overline{G}$ 
plays an essential role in the relationship 
between $G$ and $\overline{G}$ (see Section 3).   

We remind the reader that $\overline{G}$ is of 
{\em countable cofinality} if $\overline{G}$ can be 
presented as the union of an increasing $\omega$-chain 
of proper subgroups. 
This property frequently arises in ivestigations 
of automorphsm (homeomorphism) groups of countable 
structures (Polish spaces).  
We recommend the reader the paper \cite{sst} 
(and the recent paper of Ch.Rosendal \cite{rosendal}) 
for a small survey of the subject.  

In Section 4 we study the question if $\overline{G}$ 
can be $\omega$-categorical.  
Here we concentrate on the case when $G$ is soluble. 
Proposition \ref{extensions} partially explains to us 
why this case is so important. 
Using the ideas of \cite{archmac} we show that 
the condition that $\overline{G}$ is $\omega$-categorical, 
is very restrictive. 
In particular it implies that when $G$ is soluble 
and residually finite, the group $\overline{G}$ 
is nilpotent-by-finite. 
Coverings in the lattice of definable subgroup play 
an essential role both in these respects and in 
the case of cofinality.   

In Section 5 we study measurable subsets of $\overline{G}$ 
with respect to a corresponding Haar measure.
We characterize the case when the set of pairs of 
commuting elements of $\overline{G}$ has a fixed 
positive measure in any compact subgroup of $\overline{G}$.  
We also describe some aproximating function naturally 
arising in these respects. 

It looks likely that aspects of measure and 
category can be exploited much further. 
In particular we hope that our approach 
can be used in some well-known problems 
on $\omega$-categorical groups (BCM-conjecture 
and Apps-Wilson conjecture). 

\bigskip

Our basic algebraic and model-theoretic notation 
is standard (the same as in \cite{apps}).
For example $[x,y]$ is an abbreviation for $x^{-1}y^{-1}xy$ 
and $G^{(0)}=G$, $G^{(1)}=[G,G]$,...,$G^{(i+1)}=[G^{(i)},G^{(i)}]$,... 
.
When $\phi (\bar{x},\bar{y})$ is a formula, $M$ is a structure and 
$\bar{b}$ is a tuple from $M$, then we denote by $\phi (M,\bar{b})$ 
the set $\{ \bar{a}\in M : M\models \phi (\bar{a},\bar{b})\}$. 
\parskip0pt 

It is assumed in Section 4 that the reader is already 
acquainted with the most basic notions and facts of 
nilpotent groups and representations of finite groups 
(for example, commutator identities, Fitting subgroups
\footnote{maximal normal nilpotent subgroup}
, Maschke's theorem and Schur's lemma). 
  
We will also use some basic theory of profinite groups. 
Among general topological arguments (as in \cite{WilsonProf}, 
Section 0.3), we will use profinite Sylow theory. 
We will permanently consider the situation when 
a Polish group $G$ acts continuously on a Polish space $X$. 
Then we assume that the reader understands some obvious 
consequences of this situation. 
For example, the stabilizer of an element of the space  
(in particular the centralizer of an element in 
a topological group) is closed. 
If $G$ is compact then each $G$-orbit in $X$ 
is compact and closed. \parskip0pt 

In fact our model-theoretic arguments are restricted 
to some basic properties of $\omega$-categorical 
theories (outside the title of the paper we prefer 
$\omega$-categorical to $\aleph_0$-categorical). 
We now give some preliminaries concerning 
$\omega$-categorical groups. 
Let $G$ be $\omega$-categorical. 
Since all $Aut(G)$-invariant relations of $G$ 
are definable, all characteristic subgroups 
of $G$ are definable. 
It is also worth noting that for any natural 
$k$ the group $G$ has only finitely many 
$Aut(G)$-invariant $k$-ary relations.
In particular $G$ is {\em uniformly 
locally finite}: there is a function 
$f:\omega \rightarrow \omega$ such that every 
$n$-generated subgroup of $G$ is of size $\le f(n)$. 
\parskip0pt 

By $\langle g\rangle ^{G}$ we denote 
the normal subgroup of $G$ generated by $g$.  
We mention Lemma 3.2 from \cite{m1}.

\begin{lem}
If $G$ is an $\omega$-categorical group, then there is 
a formula $\psi (x,y)$ such that $G\models \psi(g,h)$
if and only if $h\in \langle g\rangle ^{G}$.
\end{lem}

This is a consequence of the fact that $G$ has 
{\em finite conjugate spread} (see \cite{apps}): 
there exists a number $m$ such that for any 
$g\in G$ any element of $\langle g\rangle^G$ is 
a product of $\le m$ conjugates of $g^{\varepsilon}$, 
$\varepsilon\in\{ -1,1\}$. 
\parskip0pt 

Theorem A of \cite{apps} will play an important 
role in our arguments. 
It states that an $\omega$-categorical group $G$ 
has a finite series $1=G_0 <G_1 <...<G_n =G$ with 
each $G_i$ characteristic in $G$, and with each 
$G_i /G_{i-1}$ either elementary abelian, or 
isomorphic to some Boolean power $P^B$ with finite 
simple non-abelian $P$ and an atomless Boolean ring $B$, 
or an $\omega$-categorical characteristically 
simple non-abelian $p$-group for some $p$. 
It is still unknown if the third case can happen 
(the Apps-Wilson conjecture states that it cannot happen). 
\parskip0pt 
 
Let us consider the case of Boolean powers in more 
detail (following \cite{apps} and \cite{macros}). 
Let $B$ be a Boolean ring (that is, a commutative 
ring in which every element is an idempotent). 
Then $B$ can be identified with the ring of compact open 
subsets of the Stone space $X=S(B)$ of $B$. 
We remind the reader that the space $S(B)$ consists 
of all maximal ideals of $B$, and the ring of all compact 
open subsets of $S(B)$ is considered with respect to 
symmetric difference + and intersection $\cdot$. 

The {\em Boolean power} of $G$ by $B$, denoted by $G^B$, 
is defined to be the set of functions $X\rightarrow G$ 
which are continuous with compact support ($G$ is disctete).  
When $g\in G$ and $A\in B$, then by $g_A$ we denote the 
element of $G^B$ with support $A$ which takes each 
$x\in A$ to $g$. 
Then every element of $G^B$ can be written uniquely 
(up to order) in the form 
$(g_1 )_{A_1}\cdot ...\cdot (g_n )_{A_n}$ with 
$g_1 ,...,g_n$ distinct elements of $G\setminus\{ 1\}$; 
and $A_1 ,...,A_n$ pairwise disjoint elements of 
$B\setminus \{ 0\}$. 
A theorem of Waszkiewicz and W\c{e}g{\l}orz claims that 
when $G$ is $\omega$-categorical and $B$ has finitely 
many atoms, the group $G^B$ is $\omega$-categorical. 
\parskip0pt 

It is worth noting that there are exactly two 
countable atomless Boolean rings: $V_0$ without 
an identity element and $V_1$ with an identity element. 
These are $\omega$-categorical, and any 
$\omega$-categorical countable Boolean ring is 
isomorphic to a product of one of these with 
a finite field of sets. 
   
We will concentrate on the case when $G$ is 
a finite simple non-abelian group. 
In this case any normal subgroup of $G^B$ 
is of the form $G^S$ where $S$ is an ideal 
of $B$ (see \cite{apps2}, Proposition 5.1).  
Theorem 2.3 of \cite{apps} states that $G^B$  
does not have proper characteristic subgroups 
if $B$ is finite or a countable atomless Boolean 
ring (i.e. $R=V_0$ or $R=V_1$). \parskip0pt   

One of the statements of Theorem B of \cite{apps} 
is that if a countable characteristically simple 
non-abelian group has a subgroup of finite index, 
then it is isomorphic to some $P^B$ where $P$ is 
a finite simple non-abelian group and $B$ is 
a countable characteristically simple Boolean ring.  

We will also use Boolean powers of rings.  
If $R$ is such a ring then $R^B$ is defined by 
the same definition as in the case of groups. 
In fact we will use {\em filtered Boolean powers} 
of rings. 
We remind the reader that a filtered Boolean power 
of $R$ (of a group $G$ resp.) with respect to 
a Boolean space $X=S(B)$, where $B$ is 
a Boolean algebra, is defined by a finite 
sublattice $L$ of closed subsets of $X$ and 
an isomorphism $\tau$ from $L$ to a sublattice 
of subrings of $R$ (subgroups of $G$ resp.). 
The corresponding filtered Boolean power is 
the substructure of $R^B$ (of $G^B$ resp.) 
consisting of those $f\in R^B$ such that for 
each $C \in L$, $f(C )\subseteq \tau (C )$. 
Denote this structure by $C(X,R,\tau )$ 
(by $C(X,G,\tau )$ resp.).  
Note that each closed subsets $C$ from $L$ 
naturally corresponds to some ideal $I_C$ of 
$B$ with $C = \{ I\in X: I_C \subseteq I\}$. 
Then the structure $( B, I_C )_{C \in L}$ is 
called an {\em augmented Boolean algebra}.   

The following Macintyre-Rosenstein description 
of $\omega$-categorical rings (see \cite{macros}) 
is one of the basic results in the area.   
\begin{quote} 
Let $R$ be an $\omega$-categorical 
ring with 1 with no nilpotent elements. 
Then $R$ is a direct product of finitely many 
0-definable rings $R_i$ and each $R_i$ is isomorphic 
to a filtered Boolean power of a finite field 
$F_i$ with respect to an appropriate Boolean 
space $X_i$ and a map $\tau_i$ to a lattice of 
subfields of $F_i$. 
The corresponding augmented Boolean algebras $B_i$ 
are interpretable in $R$ (and thus $\omega$-categorical). 
\end{quote} 

When $R=F^B$ is a Boolean power of a finite field, 
all ideals of $R$ can be described in a very similar 
way as in the case of a Boolean power of a finite simple 
non-abelian group (Proposition 5.1 of \cite{apps2}): 
any ideal of $F^B$ is of the form $F^S$ where $S$ 
is an ideal of $B$. 
We think that this is a folklour fact. 
On the other hand the proof of Proposition 5.1 from 
\cite{apps2} works for $F^B$ with obvious changes 
(for example, commutators $[g_i ,h] \not=1$ should be 
replaced by products $g_i h \not=0$). 

\section{Topological completion}

A subgroup $H$ of an infinite group $G$ is called {\em inert in} $G$ 
if for every $g\in G$, $H\cap gHg^{-1}$ is of finite index in $H$. 

\begin{prop} \label{Belyaev} 
Let $G$ be a countably infinite uniformly locally finite (for example, 
countably categorical) group. \\ 
(1) Then for every finite subgroup $F<G$ there is a finite nontrivial 
subgroup $K$ such that $F\le N_{G}(K)$ and $K\cap F =1$. \\ 
(2) Every finite subgroup of $G$ is contained in an infinite 
resudually finite inert subgroup of $G$. 
\end{prop} 

{\em Proof.} 
We use Corollary 1.7 of \cite{belyaev}:  
{\em A locally finite uncountable group satisfies statement (1) of the 
formulation}. 
If $G$ contains a finite subgroup $F$ which does not satisfy (1), 
then $F$ does not satisfy (1) in any uncountable elementary extension 
$G_1$ of $G$. 
Since $G$ is uniformly locally finite, $G_1$ is locally finite. 
Thus we have a contradiction. 

Statement (2) follows fom (1) by Theorem 1.2 of \cite{belyaev} 
$\Box$ 

\bigskip  

Let $H$ be an inert residually finite subgroup of $G$ and let 
$$
X=\{ gK : K \mbox{ is commensurable with }H \mbox{ and }g\in G\}.
$$ 
Since $G$ is countable, $X$ is countable too. 
By Theorem 7.1 of \cite{belyaev} the action of $G$ 
on $X$ by left multiplication defines an embedding 
of $G$ into $Sym(X)$ such that the closure of $G$ 
in $Sym(X)$ is a locally compact subgroup $\overline{G}$ 
and the closure of any $K$ commensurable with $H$ is 
a compact subgroup of $\overline{G}$. 
We call the group $\overline{G}$ the {\em completion} 
of $G$ with respect to the $G$-space $(G,X)$ 
(or with respect to $H$). \parskip0pt 

It is easy to see that when $G$ is residually finite, 
the completion $\overline{G}$ (with respect to $G$) 
is the profinite completion of $G$. 
This observation can be developed as follows. 
Let us enumerate $X=\{ \xi_1 ,\xi_2 ,...,\xi_n ,...\}$ 
and present $X$ as an increasing sequence of 
$G$-invariant subsets $X_n =G\xi_1 \cup ...\cup G\xi_n$, 
$n\in \omega$.   
Let $G_n$ be the permutation group induced by $G$ on $X_n$.  
Here we consider $G_n$ as a topological group 
with respect to the topology of the action on $X_n$.  
The following lemma will be helpful below. 

\begin{lem} \label{sbgrps}
Assume that $G$ is residually finite and the subgroup 
$H$ of the construction coincides with $G$. \\ 
(1) The group $\overline{G}$ is the inverse limit of 
the system of finite permutation groups $G_n$ with respect to 
the system of appropriate projections (induced by restrictions). \\ 
(2) Let $N_1 <N_2 <G$ be subgroups of $G$ with $N_1 \lhd G$. 
Then the inverse system 
$\{...\leftarrow (G_n \cap N_2)/(G_n \cap N_1 )\leftarrow ...\}$ 
induced by 
$\{...\leftarrow G_n \leftarrow G_{n+1} \leftarrow ...:n\in\omega\}$ 
has the limit isomorphic to the quotient 
$\overline{N}_2 /\overline{N}_1$ of the corresponding 
closures in $\overline{G}$. 
\end{lem} 

{\em Proof.} Since each $\xi_n$ is a coset of a subgroup 
of finite index, statement (1) is obvious. 
Now the fact that in (2) the group $\overline{N}_2$ 
($\overline{N}_1$ resp.) is the inverse limit of the 
system $\{...\leftarrow G_n \cap N_2 \leftarrow ...\}$ 
(or $\{...\leftarrow G_n \cap N_1 \leftarrow ...\}$ resp.) 
is standard. 
Thus statement (2) of the lemma becomes standard too. 
$\Box$ 

\bigskip 

We now collect some basic facts about our construction. 

\begin{prop} \label{Aext}
Let $G$ be an $\omega$-categorical group. 
Any finitely generated subgroup of the closure $\overline{G}$ 
is a homomorphic image of some finite subgroup of $G$. 
In particular, $\overline{G}$ is locally finite. 
If $G$ is $k$-step nilpotent (soluble), then $\overline{G}$ is 
$k$-step nilpotent (soluble) too. 
\end{prop} 

{\em Proof.} 
We prove a slightly stronger statement. 
Let $\{ c_1 ,...,c_k, a_{1},...,a_l \}$ be a finite 
subset of $\overline{G}$, where $c_1 ,...,c_k \in G$. 
We consider $\bar{c}\bar{a}$ as a tuple of permutations on $X$. 
Let $(a_{i,j})_{j\in\omega}$ be a sequence of elements of $G$ 
converging to $a_i$ in $Sym(X)$ where $1\le i\le l$. 
Since $G$ is uniformly locally finite we may assume that 
for all $s\not= t$ the map $a_{i,s}\rightarrow a_{i,t}$, $1\le i\le l$, 
extends to a $\bar{c}$-stabilizing isomorphism between the groups 
$\langle \bar{c},a_{1,s},...,a_{l,s}\rangle$ and 
$\langle \bar{c},a_{1,t},...,a_{l,t}\rangle$. 
This obviously implies that $\langle \bar{c},\bar{a}\rangle$ is 
a homorphic image of $\langle \bar{c}, a_{1,s},...,a_{l,s}\rangle$. 
The rest is obvious. 
$\Box$ 

\bigskip 

We will also use the following lemma. 

\begin{lem} \label{com} 
If in the situation above the group $\overline{G}$ 
is compact, then for every characteristic subgroups 
$K,L<G$ the commutator subgroup $[\overline{K},\overline{L}]$ 
is the closure of $[K,L]$ in $\overline{G}$.   
\end{lem} 

{\em Proof.} 
Let $g_1 ,...,g_k \in \overline{K}$, 
$g'_1 ,...,g'_k \in \overline{L}$ and 
$w(g_1 ,...,g_k ,g'_1 ,...,g'_k )$ be a word cosisting 
of  commutators $[g_i ,g'_i ]$.  
Let $g_{ij}$ be a sequence from $K$ converging to 
$g_i$, $i\le k$, and $g'_{ij}$ be a sequence from $L$ 
converging to $g'_i$, $i\le k$.
Then $w(g_{1j} ,...,g_{kj} ,g'_{1j},...,g'_{kj} )$, 
$j\rightarrow \infty$,  converges to 
$w(g_1 ,...,g_k ,g'_1 ,...,g'_k )$. 
Thus $w(g_1 ,...,g_k ,g'_1 ,...,g'_k )$ belongs to 
the closure of $[K,L]$. 

To see the converse take a sequence $g_j \in [K,L]$ 
converging to some $g\in \overline{G}$. 
Since $G$ is $\omega$-categorical, there is $l\in \omega$ 
such that each $g_j$ can be written as a product of $l$ 
commutators 
$[h_{1,j},h'_{1,j}]^{\varepsilon_1}\cdot ...\cdot [h_{l,j},h'_{l,j}]^{\varepsilon_l}$ 
with $h_{i,j}\in K$, $h'_{i,j}\in L$ and $\varepsilon_i \in\{ -1,1\}$. 
Since $\overline{G}$ is compact we may assume that the 
sequence of tuples $\bar{h}_j \bar{h}'_j$ is convergent 
to some $h_1 ,...,h_l ,h'_1 ,...,h'_l$ from $K\cup L$. 
Thus 
$g = [h_1 ,h'_1 ]^{\varepsilon_1}\cdot ...\cdot [h_l ,h'_l ]^{\varepsilon_l}$.    
$\Box$ 

\bigskip 

The following lemma shows that in the situation of 
Lemma \ref{com} the group $\overline{G}$ also has 
finite conjugate spread.   

\begin{lem} \label{spread} 
Assume that $G$ is $\omega$-categorical and 
residually finite (thus $\overline{G}$ is profinite). 
There is a number $m$ such that for any $g\in\overline{G}$ 
any element of $\langle g\rangle ^{\overline{G}}$ is 
a product of $\le m$ conjugates of $g^{\varepsilon}$, 
$\varepsilon\in\{ -1,1\}$.  
\end{lem} 

{\em Proof.} 
Since $G$ is $\omega$-categorical, there is a number 
$m$ such that for any $g\in G$ any element of 
$\langle g\rangle ^{G}$ is a product of $\le m$ 
conjugates of $g^{\varepsilon}$, $\varepsilon\in\{ -1,1\}$.  
If $g\in \overline{G}$ and 
$g' =\prod_{i\le l} h^{-1}_i g^{\varepsilon_i}h_i$, 
then find sequences $g_j \in G$, $j\in \omega$, and 
$h_{i,j}\in G$, $i\le l$, $j\in\omega$, converging 
to $g$ and $h_i$ respectively. 
Then each 
$g'_j = \prod_{i\le l} h^{-1}_{ij} g^{\varepsilon_i}_j h_{ij}$ 
is a product of $m$ conjugates of $g^{\varepsilon}_j$: 
$g''_{1j}$, ...,$g''_{mj}$ from $G$.  
Since $\overline{G}$ is compact we may assume 
(taking a subsequence if necessary) that 
$g''_{1j}$, ...,$g''_{mj}$ $j\in \omega$, converges 
to some conjugates of $g^{\varepsilon}$: 
$g''_{1}$, ...,$g''_{m}$ from $\overline{G}$. 
Thus $g' = g''_{1}\cdot ...\cdot g''_{m}$.   
$\Box$

\section{Cofinality} 

Let $G$ be an $\omega$-categorical group and $\overline{G}$ 
be the completion of $G$ with respect to an appropriate $H$ 
as it was defined in Section 2. 
In this section we investigate cofinality of $\overline{G}$. 

\begin{prop} \label{nncomp}  
Under circumstances above the following statements hold. \\ 
(1) If the group $\overline{G}$ is not compact, then $\overline{G}$ 
has countable cofinality in the following stronger sense: 
$\overline{G}$ is the union of an $\omega$-chain of proper open 
subgroups ($\overline{G}$ has countable topological cofinality). \\  
(2) If the group $G$ is locally soluble, then $\overline{G}$ 
has countable cofinality. 
\end{prop} 

{\em Proof.} 
Let $H$ be an inert residually finite subgroup of $G$, 
$$ 
X=\{ gK : K \mbox{ is commensurable with } H \mbox{ and } g\in G\},
$$
and $\overline{G}$ be the corresponding closure. 
 
(1) 
It suffices to prove that $\overline{G}$ is not compactly generated 
(see p.9 of \cite{rosendal}). 
Let $F\subset \overline{G}$ be compact and $gK$ be any element of $X$.  
Since $F$ is compact and the topology of $\overline{G}$ is defined 
by the topology of $Sym(X)$, the set $\{ fgK : f\in F\}$ 
(of images of $gK$ under $F$) is finite (say $f_1 gK ,...,f_l gK$). 
Thus $Fg\subset \bigcup f_i gK$ and $F$ is contaned in the subgroup 
$K_1 =\langle f_1 ,...,f_k ,g, K\rangle$. 
\parskip0pt 

Let $K_0 = \bigcap \{ K^h : h\in \langle f_1 ,...,f_k , g\rangle\}$. 
Then $K_0$ is $\langle \overline{f},g\rangle$-normalized and 
of finite index in $K$. 
Let $h_1 ,...,h_l$ represent all nontrivial cosets of $K/K_0$.  
Since $G$ is locally finite, the group $K_0$ is of finite 
index in $\langle \overline{f},g,\overline{h},K_0 \rangle$. 
We see that $K$ is of finite index in 
$K_1 =\langle f_1 ,...,f_k ,g, K\rangle$.  
Thus $K_1 \in X$ and by Theorem 7.1 of \cite{belyaev}
the closure $\overline{K}_1$ is compact. 
We now see that $F$ cannot generate $\overline{G}$. 

(2) By statement (1) we may assume that $G$ is 
residually finite and $\overline{G}$ is compact. 
In this case $\overline{G}$ is prosoluble and 
$[\overline{G},\overline{G}]$ is a proper closed 
subgroup of $\overline{G}$ (Lemma \ref{com}). 
We claim that $\overline{G}$ is really soluble. 

To see this we apply Theorem A of \cite{apps}. 
By this theorem $G$ has a finite series 
$1=G_0 <G_1 <...<G_n =G$ with each $G_i$ characteristic in $G$, 
and with each $G_i /G_{i-1}$ either elementary abelian, or 
isomorphic to some Boolean power $P^B$ with finite simple 
non-abelian $P$ and atomless Boolean ring $B$, 
or an $\omega$-categorical characteristically simple 
non-abelian $p$-group for some $p$. 
By Theorem B of \cite{apps} an infinite $\omega$-categorical 
characteristically simple non-abelian $p$-group is not 
residually finite. 
Thus this case cannot be realized as $G_i /G_{i-1}$.  
In the case when $G_i /G_{i-1}$ is of the form $P^B$ 
the group $P$ is a homomorphisc image of a subgroup 
of $G$. 
Since $G$ is locally finite and locally soluble, 
the group $P$ cannot be simple non-abelian. 
As a result we see that all factors 
$G_i /G_{i-1}$ are abelian and the grpoup $G$ is soluble. 
By Proposition \ref{Aext} the group $\overline{G}$ is soluble too. 
\parskip0pt 

Consider the chain of derived subgroups 
$\overline{G}>[\overline{G},\overline{G}]>[[\overline{G},\overline{G}],[\overline{G},\overline{G}]]>...$.  
Since $\overline{G}$ is soluble there is $n$ such that the 
subgroup $\overline{G}^{(n)}$ is of finite index in $\overline{G}$, 
but $\overline{G}^{(n+1)}$ is of infinite index in $\overline{G}^{n}$. 
In this case the cofinality of $\overline{G}$ coincides with 
the cofinality of $\overline{G}^{(n)}$. 
Since an infinite abelian group has a countably infinite 
quotient group (Lemma 4.6 of \cite{rosendal}) there is 
a normal subgroup $L \le\overline{G}^{(n)}$ with countable 
$\overline{G}^{(n)}/L$ such that $\overline{G}^{(n+1)}\le L$.  
Let $g_1 ,...,g_n ,...$ represent all cosets of $\overline{G}^{(n)}/L$. 
Since $\overline{G}^{(n)}$ is locally finite, the sequence 
$$
\L\le \langle L,g_1 \rangle \le ...\le \langle L,g_1 ,...,g_n \rangle \le...  
$$ 
contains a strictly increasing chain such that its 
union equals $\overline{G}^{(n)}$. 
$\Box$ 

\bigskip 

Assume that $G$ is residually finite. 
Consider the lattice of all characteristic subgroups of $G$. 
As we have already mentioned (in the proof above), by Theorems A 
and B of \cite{apps} if $N_1 <N_2$ is a {\em covering pair} from this 
lattice (i.e. there is no intermediate characteristic subgroups 
between $N_1$ and $N_2$), then $N_2 /N_1$ is either elementary abelian, 
or isomorphic to some Boolean power $P^B$ with finite simple 
non-abelian $P$ and atomless Boolean ring $B$. 
When the former case appears, the group $N_2 /N_1$ becomes 
a $G/N_2$-module under the conjugacy action.  
Note that since $N_2 /N_1$ is abelian, this action 
does not depend on representatives of $N_2$-cosets.  
The same property holds for the topological 
$\overline{G}/\overline{N}_2$-module 
$\overline{N}_2 /\overline{N}_1$ (by Lemma \ref{sbgrps}). 
\parskip0pt 

Under the circumstances as above we will say that 
$N_2 /N_1$ has {\em countable $G/N_2$-cofinality} 
if there are a subgroup $H<G$ and an 
$\omega$-chain $N_1 <U_0 <U_1 <U_2 <...$ of subgroups 
of $N_2$ such that $G=HN_2$, 
$U_0 = HN_1  \cap N_2$, $\bigcup U_i =N_2$ 
and $(U_0 /N_1 ) < (U_1 /N_1 )<(U_2 /N_1 )<... $ is 
a proper $\omega$-chain of $G/N_2$-submodules of 
$N_2 /N_1$.  
The same definition (under the replacement 
$N_2 /N_1 \rightarrow \overline{N}_2 /\overline{N}_1$) 
works in the case of $\overline{G}$.  
\parskip0pt 

This definition is slightly stronger than the natural 
version of the condition that $N_2 /N_1$ has cofinality $\omega$. 
Indeed, our definition guarantees that the sequence 
of subgroups (use the fact that each $U_i$ is an $H$-module) 
$HU_0 <HU_1 <HU_2 < ...$ is proper and the corresponding 
union coincides with $G$ (i.e. $G$ has countable cofinality). 
To see the former condition we easily reduce the problem 
to the situation when some $u\in U_{k+1} \setminus U_k$ 
is the product of some $h\in H$ and $u'\in U_k$. 
In this case $h\in N_2$ and our definition implies 
$h\in U_0\subset U_k$, a contradiction. \parskip0pt   

By Proposition 3.10 of \cite{rosendal} any compact 
Polish group has uncountable topological cofinality. 
In the following theorem we apply the term which we have 
just introduced to characterize when $\overline{G}$ has 
countable cofinality.   

\begin{thm} \label{characterisation}  
Let $G$ be a countably categorical group. 
Then the completion $\overline{G}$ has uncountable 
cofinality if and only if $G$ is residually finite 
and the following properties hold: \\ 
(1) in any series $1=G_0 <G_1 <...<G_n =G$ of 
characteristic subgroups in $G$, with each $G_i /G_{i-1}$ 
characteristically simple, the infinite quotient 
$G_m /G_{m-1}$ with maximal $m$ is isomorphic 
to some Boolean power $P^B$ with finite simple 
non-abelian $P$; \\ 
(2) for any abelian cover $N_1 <N_2$ in 
lattice of characteristic subgroups of $G$ 
if $N_2$ is of infinite index in $G$, then 
the $\overline{G}/\overline{N}_2$-module 
$\overline{N}_2 /\overline{N}_1$ 
does not have countable $\overline{G}/\overline{N}_2$-cofinality. 
\end{thm} 

We start with some preliminary lemma. 

\begin{lem} \label{prod}  
Let $P$ be a finite simple non-abelian group, 
$V_1$ ($V_0$ resp.) be the unique countable atomless 
Boolean ring with (without) identity 
and $G=P^{V_1}$ ($G=P^{V_0}$ resp.). 
Then the profinite completion $\overline{G}$ is 
isomorphic to the direct power $P^{\omega}$.    
\end{lem} 

{\em Proof.} 
We apply Proposition 5.1 of \cite{apps2}, which describes 
normal subgroups of Boolean powers. 
According to this result the normal subgroups of $G$ 
are just those of the form $P^S$, for $S$ an ideal of 
$B= V_1$ ($=V_0$ resp.). 
It is easy to see that this subgroup is the kernel of 
the homomorphism $P^B \rightarrow P^{B/S}$ mapping 
each $(g_1 )_{A_1}\cdot \cdot \cdot (g_k )_{A_k}$  
(with $A_1 ,...,A_k$ disjoint) to 
$(g_1 )_{A_1 +S}\cdot \cdot \cdot (g_k )_{A_k +S}$.  
When $P^{B/S}$ is finite, it is naturally isomorphic  
to the group $P^m$ where $m$ is the number of atoms of $B/S$. 
Now it is easy to see that the corresponding inverse system 
consists of projections between groups of the form $P^m$. 
The corresponding limit is the direct power $P^{\omega}$. 
$\Box$ 

\bigskip  

We will also use the following theorem 
of S.Koppelberg and J.Tits \cite{KoTits}: 
\begin{quote}
Let $F$ be a non-trivial finite group. 
Then the power $F^{\omega}$ has uncountable cofinality 
if and only if $F$ is perfect. 
\end{quote} 

{\em Proof of Theorem \ref{characterisation}.} 
Necessity follows from Theorem B of \cite{apps} and 
Proposition \ref{nncomp}. 
Indeed assuming the contrary we are left to the case when  
$G$ is residually finite (by Proposition \ref{nncomp}(1)) and 
condition (2) is satisfied (see the discussion above 
concerning countable $G/N_2$-cofinality).  
Theorem B of \cite{apps} guarantees that $G_m /G_{m-1}$ 
is not an infinite non-abelian characterically simple $p$-group. 
Thus, when (1) does not hold, $\overline{G}$ has a subgroup 
of finite index with an infinite abelian quotient. 
Now Proposition \ref{nncomp}(2) implies that this quotient 
and $\overline{G}$ have countable cofinality.  
\parskip0pt 

To prove sufficiency consider $1=G_0 <G_1 <...<G_n =G$, 
any series of characteristic subgroups in $G$, with each 
$G_i /G_{i-1}$ characteristically simple. 
Find $m\le n$ such that $|G:G_m |$ is finite (say $k$) 
and $|G_m :G_{m-1}|$ is infinite. 
By our assumptions the infinite quotient $G_m /G_{m-1}$ is 
isomorphic to some Boolean power $P^B$ with finite simple 
non-abelian $P$ and $B=V_0$ or $=V_1$. 
\parskip0pt 

Consider all finite quotients of $G$. 
If $L$ is a finite quotient of $G$ with $|G:G_m|<|L|$, then 
the index of the image of $G_m$ in $L$ is also bounded by $k$. 
Thus the closure $\overline{G}_{m}$ in $\overline{G}$ 
is of index $k$ in $\overline{G}$ and the closure 
$\overline{G}_{m-1}$ is a normal subgroup such that 
$\overline{G}/\overline{G}_{m-1}$ is a profinite completion 
of $G/G_{m-1}$ (by Lemma \ref{sbgrps}). 
By Lemma \ref{prod} the quotient $\overline{G}_m /\overline{G}_{m-1}$ 
is of the form $P^{\omega}$ for $P$ as above. 

We now assume that $\overline{G} =\bigcup H_i$ 
for some $\omega$-chain $1 <H_1 <...<H_i <...$ 
of subgroups of $\overline{G}$.   
Choose the minimal $l\le m$ with the property 
that for some $i\in \omega$, $H_i \overline{G}_{l}=\overline{G}$.  
It is clear that we may assume that $l>0$. 
Consider the sequence 
$(H_i \overline{G}_{l-1})/\overline{G}_{l-1} \cap \overline{G}_l /\overline{G}_{l-1}$, 
$i\in \omega$ . 
We may assume this sequence is strictly increasing and 
$H_1 \overline{G}_l =\overline{G}$. \parskip0pt 

If $\overline{G}_l /\overline{G}_{l-1}$ is abelian, then $l<m$, 
i.e. $\overline{G}/\overline{G}_l$ is infinite. 
Representing $\overline{G}/\overline{G}_l$ as 
$H_1 \overline{G}_l /\overline{G}_l$ we easily see that each 
$(H_i \overline{G}_{l-1} )/\overline{G}_{l-1}\cap \overline{G}_{l}/\overline{G}_{l-1}$ 
is a $\overline{G}/\overline{G}_l$-module. 
Taking $U_i = H_i \overline{G}_{l-1}\cap \overline{G}_l$ 
we obtain a contradiction with property (2) of the formulation. 
\parskip0pt 

If $\overline{G}_l /\overline{G}_{l-1}$ is not abelian then we see as 
above that $\overline{G}_l /\overline{G}_{l-1}$ is isomorphic to some 
product of the form $P^{\omega}$, where $P$ is finite, simple 
and non-abelian. 
On the other hand the sequence 
$(H_i \overline{G}_{l-1} )/\overline{G}_{l-1}\cap \overline{G}_{l}/\overline{G}_{l-1}$, 
$i\in\omega$, covers $\overline{G}_{l}/\overline{G}_{l-1}$, 
i.e. $\overline{G}_{l}/\overline{G}_{l-1}$ has countable cofinality. 
This contradicts a theorem of Koppelberg and Tits (see \cite{KoTits}).   
$\Box$ 
\bigskip 

The methods applied above can also work for some other 
questions naturally arising in our approach. 
Assume that $G$ is an $\omega$-categorical, residually finite 
group. 
As we already know, for any normal subgroups $N_1 <N_2 <G$ 
the quotient $\overline{N}_2 /\overline{N}_1$ of their 
closures in the completion $\overline{G}$ is the profinite 
completion of $N_2 /N_1$.  
As $G$ is residually finite, the lattice of 
characteristic subgroups of $G$ does not contain 
any covers $N_1 <N_2$ with $N_2 /N_1$ non-abelian and  
characteristically simple $p$-group (Theorem B of \cite{apps}). 
What are possibilities for $\overline{N}_2 /\overline{N}_1$ ?

\begin{cor} \label{smallness}
Assume that $G$ is an $\omega$-categorical, residually finite 
group. 
If for every cover $N_1 <N_2$ of the lattice of 
characteristic subgroup of $G$, the completion of $N_2 /N_1$ 
can be presented as an inverse limit which is small 
in the sense of Newelski, then $\overline{G}$ is soluble. 
\end{cor} 

Let us recall necessary preliminaries from \cite{newelski}. 
Let $\Gamma$ be a profinite group and 
$G_0 \leftarrow G_1 \leftarrow ...$ be 
the corresponding inverse system. 
Any automorphism of $\Gamma$ fixing all $G_i$ (as sets)
is called a {\em profinite automorphism} of $\Gamma$. 
By $Aut^{*}(\Gamma )$ we denote the group of all profinite 
automorphisms of $\Gamma$. 
L.Newelski defines that $\Gamma$ is {\em small} if 
the number of all $Aut^* (\Gamma )$-orbits on $\Gamma$ is 
at most countable. 
The main conjecture in these respects states that small 
profinite groups have open abelian subgroups. 

\bigskip 

{\em Proof of Corollary \ref{smallness}} 
To prove the corollary it is enough to prove that Boolean 
products $P^B$ with $P$ finite, non-abelian and simple, 
cannot be realized by quotients of covers of the lattice 
(then by $\omega$-categoricity we see that the lattice 
is finite and by Theorem A of \cite{apps} it has only 
abelian covers, thus $G$ is soluble).  
\parskip0pt 

Assume that $N_1 <N_2$ is a cover such that 
$N_2 /N_1 \cong P^B$, where $P$ is as above. 
Then by Theorem A of \cite{apps} $B$ is an atomless 
Boolean ring, i.e. $V_0$ or $V_1$ in our notation above. 
By Lemma \ref{prod} the profinite completion of 
$N_2 /N_1$ is isomorphic to the product $P^{\omega}$. 
By Remark 4.2 of \cite{newelski}, the group 
$P^B$ is not small. 
$\Box$

\section{Categoricity of the completion} 

Let $G$ be an $\omega$-categorical group. 
We will assume here that $G$ is soluble. 
Then by the main theorem of \cite{archmac} $G$ is nilpotent-by-finite 
or interprets the countable atomless Boolean algebra. 
In the former case $\overline{G}$ is nilpotent-by-finite too. 
In this section we show that the latter case is impossible if 
$G$ is residually finite and $Th(\overline{G})$ is $\omega$-categorical.   

We start our investigations with an observation that 
tricks of the previous sections are still helpful for 
the question when $Th(\overline{G})$ is $\omega$-categorical   
(in fact we will use them in the main theorem of the section).  
Consider $\omega$-categorical residually finite $G$ 
acting on $X=\{ \xi_1 ,\xi_2 ,...,\xi_n ,...\}$ under 
the action defined as in Section 2.  
We present $X$ as an increasing sequence of $G$-invariant 
subsets $X_n =G\xi_1 \cup ...\cup G\xi_n$, $n\in \omega$.   
Then $\overline{G}$ is the inverse limit of the permutation 
groups $(G_i ,X_i )$ induced on all $X_i$ by $G$ and with 
appropriate projections.  
Lemma \ref{sbgrps} describes the case when 
$N_1 <N_2 <G$ are normal subgroups of $G$. 
Then by Lemma \ref{sbgrps} the inverse system 
$\{ ...\leftarrow (G_i \cap N_2)/(G_i \cap N_1 )\leftarrow ...\}$ 
induced by $\{ ...\leftarrow (G_i ,X_i )\leftarrow ...\}$ has 
the limit isomorphic to the quotient $\overline{N}_2 /\overline{N}_1$ 
of the corresponding closures in $\overline{G}$. 

\begin{prop} \label{extensions}  
Le $G$ be an $\omega$-categorical residually finite group 
and $N_1 <N_2< G$ be a cover in the lattice of characteristic 
subgroups. 
If $Th(\overline{G})$ is $\omega$-categorical and the closures 
$\overline{N}_1$ and $\overline{N}_2$ are definable in $\overline{G}$ 
then $\overline{N}_2 /\overline{N}_1$ is abelian. 
\end{prop} 

{\em Proof.} 
If the cover $N_1 <N_2$ is not abelian, then by Theorems A and B 
of \cite{apps}, $N_2 /N_1 \equiv P^B$, where $P$ is finite, simple 
and non-abelian, and $B$ is an atomless Boolean ring, 
i.e. $V_0$ or $V_1$ in our notation above. \parskip0pt 

By Lemma \ref{sbgrps} the quotient $\overline{N}_2 /\overline{N}_1$ 
is the profinite completion of $N_2 /N_1$. 
On the other hand by Lemma \ref{prod} the profinite completion 
of $N_2 /N_1$ is isomorphic to the product $P^{\omega}$.  
By the assumptions of the proposition, 
$\overline{N}_2 /\overline{N}_1$ is definable in $\overline{G}$. 
Since $P^{\omega}$ cannot be $\omega$-categorical \cite{ros}, 
we see that $\overline{G}$ is not $\omega$-categorical. 
$\Box$ 

\bigskip 

Note that this proposition somehow suggests that 
in the $\omega$-categorical residually finite case the group 
$\overline{G}$ must be very similar to soluble groups: 
if there is a series of characteristic closed subgroups 
as in the Apps' theorem, then $\overline{G}$ is soluble. 
The main theorem of the section is devoted to this case. 

\begin{thm} \label{4main} 
Let $G$ be an $\omega$-categorical, residually finite, soluble 
group such that the profinite completion $\overline{G}$ has 
$\omega$-categorical elementary theory.   
Then both $G$ and $\overline{G}$ are nilpotent-by-finite. 
\end{thm} 

The proof of this theorem uses ideas of 
Proposition \ref{extensions} and the paper \cite{archmac}. 
We will assume $G$ is a soluble, residually finite, 
$\omega$-categorical group which is not nilpotent-by-finite. 
Then by Lemma \ref{Aext} $\overline{G}$ is soluble and is not 
nilpotent-by-finite. 
In order to obtain a contradiction assume 
that $\overline{G}$ is $\omega$-categorical. 
We start with the following lemma. 

\begin{lem} \label{lm1}  
Let $G$ be an $\omega$-categorical, soluble, residually 
finite group which is not nilpotent-by-finite.  
If the elementary theory of the profinite completion 
$\overline{G}$ is $\omega$-categorical, then there are 
closed characteristic subgroups $S_1 <S_2 <\overline{G}$ 
and $T_1 <T_2 <\overline{G}$ such that $S_2 /S_1$ 
is abelian over its centre, $[S_2 /S_1 ,S_2 /S_1 ]$ is finite, 
$T_2 /T_1$ is an $S_2 /S_1$-module under the conjugacy 
action and the corresponding semidirect product of 
$T_2 /T_1$ and $S_2 /S_1$ is not nilpotent-by-finite. 
\end{lem} 

{\em Proof.} 
In fact Archer and Macpherson have proved this lemma  
in \cite{archmac} without the condition of closedness 
of the corresponding subgroups (they did not have any 
topological assumptions).  
We analyse $\overline{G}$ by the method of \cite{archmac}. 
Find the minimal number $r$ such that the quotient 
$\overline{G}/\overline{G}^{(r)}$ is not nilpotent-by-finite, 
but $\overline{G}/\overline{G}^{(r-1)}$ so is. 
Let $A=\overline{G}^{(r-1)}$ and let $F$ be 
the Fitting subgroup of $\overline{G}/A$. 
By Lemma \ref{com} $A$ is the closure of $G^{(r-1)}$ and 
$\overline{G}^{(r)}$ is the closure of $G^{(r)}$.  
As $FA/A$ is a characteristic subgroup of 
$\overline{G}/A$ of finite index, we also have 
that $FA$ is a characteristic subgroup of 
$\overline{G}$ of finite index.   
Note that $FA$ is the closure in $\overline{G}$ of 
the preimage in $G$ of the Fitting subgroup 
(say $F_0$) of $G/G^{(r-1)}$. 
Indeed, by Lemma \ref{sbgrps} the quotient $\overline{G}/A$ 
is the profinite completion of $G/G^{(r-1)}$. 
Thus $G/G^{(r-1)}$ is nilpotent-by-finite and 
$F_0 G^{(r-1)}$ has the closure of finite index 
in $\overline{G}$.  
Since both $G/(F_0 G^{(r-1)})$ and $\overline{G}/(FA)$ 
do not have normal nilpotent subgroups, the group 
$FA$ is the closure of $F_0 G^{(r-1)}$ in $\overline{G}$
(remember that by Lemma \ref{Aext} the closure of a nilpotent 
characteristic subgroup is also nilpotent). 
\parskip0pt 

Let $F=Z_c >...>Z_0 =1$ be the lower central series of $F$. 
Then by Lemma \ref{com} each $Z_t A$ is a closed 
characteristic subgroup of $\overline{G}$. 
Find the maximal $t$ such that $(Z_t A) /\overline{G}^{(r)}$ 
has a nilpotent subgroup of finite index. 
By our assumptions $c\not= t$. 
Let $S_1\ge A$ be the preimage in $\overline{G}$ of 
the Fitting subgroup of $(Z_t A) /\overline{G}^{(r)}$. 
Thus $S_1$ is a characteristic subgroup of $\overline{G}$. 
Applying arguments as in the previous paragraph 
we see that $S_1$ is closed. \parskip0pt 

Since $S_1 /\overline{G}^{(r)}$ is nilpotent 
there is a sequence of $S_1$-invariant closed subgroups 
$\overline{G}^{(r)}=A_0 <...<A_n =A$, characteristic 
in $\overline{G}$, such that $S_1 /A_0$ centralizes 
each $A_{i+1}/A_i$.  
Here we assume that each $A_i$ is the intersection 
of $A$ with the corresponding member of the 
lower central series of $S_1 /G^{(r)}$ 
(thus we can apply Lemma \ref{com}). 
By a straightforward argument (see the induction 
step in the proof of Lemma 2.1 of \cite{archmac}) 
there is $i$ such that the semidirect product of 
$A_{i+1}/A_i$ and $(Z_{t+1}A)/S_1 )$ (under 
the conjugacy action) is not nilpotent-by-finite. 
Let $T_1 =A_i$, $T_2 =A_{i+1}$ and $S_2 = Z_{t+1} A$. 
As we already mentioned these subgroups are closed and 
characteristic in $\overline{G}$. \parskip0pt 

Note that the centralizer of $(Z_t A)/S_1$ in $(Z_{t+1}A)/S_1$ 
defines a characteristic closed subgroup $U<\overline{G}$ 
of finite index in $Z_{t+1}A$. 
Thus we may assume that $S_{2}/S_1$ centralizes 
$(Z_t A)/S_1$ (replacing $Z_{t+1} A$ by $U$ if necessary).  
$\Box$ 

\bigskip 

In the situation of the lemma consider the action of 
$S_{2}/S_1$ on $T_{2}/T_1$ in more detail. 
Since the semidirect product of $T_{2}/T_1$ and 
$S_{2}/S_1$ is not nilpotent-by-finite, the pointwise 
stabilizer of $T_{2}/T_1$ (fixator) of the conjugacy 
action of $S_{2}/S_1$ is of infinite index in the latter. 
Since $T_1$ and $T_{2}$ are characteristic 
in $\overline{G}$, the conjugacy-fixator of $T_{2}/T_1$ 
is also characteristic in $\overline{G}$. 
Since it is also closed we can enlarge $S_1$ so that 
$S_2 /S_1$ acts on $T_2/T_1$ faithfully. 
\parskip0pt

Decomposing the profinite group $S_2 /S_1$ into 
the direct sum of (finitely many) its Sylow subgroups 
we find a prime number $r$ such that the preimage in 
$\overline{G}$ of the $r$-Sylow subgroup (which is 
obviously characteristic and closed) satisfies all 
the assumptions of the lemma above for $S_2$.  
\parskip0pt 

On the other hand decomposing the profinite group 
$T_2 /T_1$ into the direct sum of (finitely many) its 
Sylow subgroups we find a prime number $p$ and 
the corresponding $p$-Sylow subgroup $P$ such that 
for the natural action of $S_2 /S_1$ on $P$ the 
corresponding semidirect product is not nilpotent-by-finite. 
Extending $T_1$ by the preimage in $\overline{G}$ of the 
complement of $P$ (which is obviously characteristic 
and closed) we still have that all the assumptions 
of the lemma above are satisfied.  
\parskip0pt 

Note that when $p=r$, the corresponding semidirect 
product is locally nilpotent, and thus by 
$\omega$-categoricity and 
\cite{wilson}, is nilpotent. \parskip0pt 

We now summarize the situation as follows. 

\begin{lem} \label{sum} 
Let $G$ be an $\omega$-categorical, soluble, residually 
finite group which is not nilpotent-by-finite. 
Let the completion $\overline{G}$ be $\omega$-categorical too.   
Then there are closed characteristic subgroups 
$S_1 <S_2 <\overline{G}$ and $T_1 <T_2 <\overline{G}$ 
such that for two prime numbers $r\not=p$, $S_2 /S_1$ 
is an $r$-group which is abelian over its center, 
$T_2 /T_1$ is an abelian $p$-group which is 
a faithful $S_2 /S_1$-module under the conjugacy action 
and the corresponding semidirect product of $T_2 /T_1$ 
and $S_2 /S_1$ is not nilpotent-by-finite. 
Moreover $[S_2 /S_1 ,S_2 /S_1 ]$ is finite. 
\end{lem} 

The following lemma is the main preliminary 
step to Theorem \ref{4main}. 

\begin{lem} \label{lm3} 
Under the assumptions of Lemma \ref{sum} there is 
a finite tuple $\bar{w}\in \overline{G}$ and there 
are closed $\bar{w}$-definable subgroups 
$S_1 <S_2 <\overline{G}$ and $T_1 <T_2 <\overline{G}$ 
such that for two prime numbers $r\not=p$, 
$S_2 /S_1$ is an infinite abelian $r$-group, 
$T_2 /T_1$ is an elementary abelian $p$-group which 
is a faithful $S_2 /S_1$-module under the conjugacy 
action. 
\end{lem} 

{\em Proof.} 
In the situation provided by Lemma \ref{sum} 
let $p^n$ be the exponent of $P=T_2 /T_1$ and 
for each $i = 0,...,n-1$ let $P_i = p^i P$ and $K_i$ be 
the centralizer of $P_i /P_{i+1}$ in $S_2 /S_1$. 
Since all $P_i T_1$ are characteristic and closed in $\overline{G}$, 
each $K_i S_1$ is characteristic and closed too. 
If all $K_i S_1$ are of finite index in $S_2$,  
then $(\bigcap K_i )$ is an $r$-group centralizing the $p$-group 
$P$ (Theorem 5.3.2 of \cite{gorenstein}) contradicting the 
assumption that the action is faithful.
Replacing $P$ by an appropriate $P_i /P_{i+1}$ and $S_2 /S_1$ 
by the corresponding $S_2 /(K_i S_1 )$ we arrange that 
$T_2 /T_1$ is a vector space over $GF(p)$. 
Denote $V=T_2 /T_1$. \parskip0pt 

We now follow the proof of Lemma 2.7 of \cite{archmac}. 
Our strategy is to show that at every step of that proof 
we build closed subgroups of $\overline{G}$. 

Assume that $S_2 /S_1$ is not abelian. 
Let $K$ be the largest finite characteristic 
subgroup of $S_2 /S_1$. 
Then according to the construction from the proof of 
Lemma \ref{lm1} we have $[S_2 /S_1, S_2 /S_1 ]\le K$. 
If the centre of $S_2 /S_1$ is infinite the proof 
is finished by replacing $S_2 /S_1$ by its centre. 
Thus passing to the centralizer of $K$ if necessary, 
we may assume that $K$ is the centre of $S_2 /S_1$. 
\parskip0pt 

By Maschke's theorem we may write $V$ as a direct 
sum of $GF(p)K$-irreducible submodules. 
Since the number of isomorphism types of such 
irreducibles is finite, by collecting them in this 
sum into isomorphism classes we may suppose that 
$V=W_1 \oplus ...\oplus W_t$, where each $W_i$ is 
a direct sum of isomorphic $GF(p)K$-irreducibles, 
and for $i\not= j$ irreducibles appearing in $W_i$ 
and $W_j$ are not isomorphic. 
Note that each $W_i T_1$ can be easily presented 
as an intersection of clopen subgroups of $\overline{G}$ 
(appearing from appropriate finite quotients of $\overline{G}$). 
By a similar argument any subgroup of the form 
$(W_{i_1}\oplus ...\oplus W_{i_l})T_1$ is closed. 

We may suppose that $V$ is a direct sum of 
isomorphic faithful irreducible $GF(p)K$-modules. 
To see this we apply the following reductions 
of Lemma 2.7 of \cite{archmac}. 
For $i= 1,...,t$ let $C_i =C_{S_2 /S_1}(W_i )$. 
It is easy to see that some group $S_2 /(C_i S_1)$ 
is infinite (using the fact that $T_2 /T_1$ is a faithful 
$S_2 /S_1$-module) and each $C_i S_1$ is closed 
in $\overline{G}$. 
It is also clear that the group 
$(\sum \{ W_i :S_2 /(C_i S_1)$ is infinite $\} )T_1$ 
is characteristic and closed in $\overline{G}$. 
We thus can identify it with $T_2$. 
If some $S_2 /(C_i S_1)$ has infinite centre we replace $V$ 
by $W_i$ and $S_2 /S_1$ by the centre of $S_2 /(C_i S_1)$. 
We thus may suppose that each $S_2 /(C_i S_1)$ 
has finite centre. 
Now replace $V$ by an appropriate $W_i$ and 
$S_2 /S_1$ by the corresponding $S_2 /(C_i S_1)$. 

Now by Schur's lemma, each non-trivial element of $K$ 
acts fixed-point-freely on $V$. 
The second claim of the proof of Lemma 2.7 from \cite{archmac} 
states that there exists a non-trivial $w\in V$ such 
that the centralizer $C_{S_2 /S_1}(w)$ is infinite. 
Since $K\setminus \{ 1\}$ acts fixed-point-freely on 
on $V\setminus \{ 0\}$, the centralizer $C_{S_2 /S_1}(w)$ 
is disjoint from $K\setminus \{ 1\}$. 
Since $S_2 /S_1$ is nilpotent of class two, 
$C_{S_2 /S_1}(w)$ is abelian. 
Since $C_{S_2 /S_1}(w)S_1$ is closed 
in $\overline{G}$ we can replace 
(adding a constant) $S_2$ by 
$C_{S_2 /S_1}(w)S_1$. 
$\Box$ 
\bigskip 

{\em Proof of Theorem \ref{4main}.} 
Let $G$ be an $\omega$-categorical, soluble, residually 
finite group which is not nilpotent-by-finite. 
Let the completion $\overline{G}$ be $\omega$-categorical too.   
As we have already known by adding finitely many constants 
we can arrange that there are closed characteristic subgroups 
$S_1 <S_2 <\overline{G}$ and $T_1 <T_2 <\overline{G}$ 
such that for two prime numbers $r\not=p$, $S_2 /S_1$ 
is an abelian $r$-group, $V= T_2 /T_1$ is 
an elementary abelian $p$-group which is a vector space 
over $GF(p)$ and a faithfull $GF(p)[S_2 /S_1 ]$-module 
under the conjugacy action. 

Furthermore by picking some additional parameters 
we may assume that 
\begin{quote} 
there is $v\in V$ such that $V=\langle v^{S_2 /S_1}\rangle$. 
\end{quote} 
Indeed, if for every $v\in V$ the centralizer 
$C_{S_2 /S_1}(v)$ is of finite 
index in the abelian $r$-group $S_2 /S_1$, 
then by Proposition 3.4 of \cite{m1} 
\footnote{if $G$ is an $\omega$-categorical group 
with a characteristic subgroup $V$ which is a vector space 
over $GF(p)$ with $G/V$ a nilpotent group with no elements 
of order $p$, and with $V$ a sum of finite-dimensional 
$GF(p)[G/V]$-modules, then the $G/V$-centralizer of $V$ 
is of finite index in $G/V$} 
(applied to the semidirect product as above) 
the group $S_2 /S_1$ cannot act faithfully on $V$. 
Thus choosing $v\in V$ with infinite vector 
space $\langle v^{S_2 /S_1}\rangle$ we see 
that this space is definable over $v$ 
(by $\omega$-categoricity) and is of the form 
$T_3 /T_1$ where $T_3 <T_2$ is closed 
(using that $\overline{G}$ is compact we see that 
the $S_2 /S_1$-orbit of $v$ is closed). 
To have a faithful action on $V$ we now replace 
$S_2 /S_1$ by its quotient over the centralizer of $v$.  

By $\omega$-categoricity there is $n\ge 1$ 
such that each element of $V$ can be written 
as a sum of at most $n$ $S_2 /S_1$-translates of $v$. 
We introduce ring operations $\oplus$ and $\otimes$ 
on $V$ as follows. 
Let $\oplus$ be the group operation on $V$. 
To define $\otimes$ let $v=\sum_{1\le i\le r} v^{h_i}$ 
and $v'=\sum_{1\le i\le s} v^{h'_i}$ be from $V$, 
where $h_i ,h'_i \in S_2 /S_1$. 
Then define $v\otimes v'$ to be 
$\sum_{1\le i\le r} \sum_{1\le j\le s} v^{h_i h'_j}$. 
It is proved in \cite{archmac} (Corollary 3.2) that 
$V$ becomes a definable ring without non-zero 
nilpotent elements.  
\parskip0pt 

Note that if $\phi$ is a homomorphism from $\overline{G}$ 
to a finite group $G_1$, then the action in $G_1$ 
induced by $\phi$ and the action of $S_2 /S_1$ on 
$T_2 /T_1$, defines as above a finite ring 
without non-zero nilpotent elements (the argument of 
Corollary 3.2 from \cite{archmac} works in this case too).  
Moreover since $\overline{G}$ is profinite the action 
of $S_2 /S_1$ on $T_2 /T_1$ is the inverse limit of 
the system of all such actions of finite groups. 
Thus the ring $V$ is the inverse limit of the system 
of corresponding finite rings. 

We now analyse the structure of $(V,\oplus ,\otimes )$ 
using the Macintyre-Rosenstein description of 
$\omega$-categorical rings from \cite{macros}.   
According this description $V$ is a direct product of 
finitely many 0-definable rings $R_i$ and each $R_i$ 
is isomorphic to a filtered Boolean power of 
a finite field $F_i$ with respect to an appropriate 
Boolean space $X_i$ and a map $\tau_i$ to a lattice 
of subfields of $F_i$ (the corresponding augmented 
Boolean algebras $B_i$  are interpretable in $V$ 
and thus $\omega$-categorical). 

Let $I$ be an ideal of $B_i$ such that $B_i /I$ 
is finite. 
Then it is easy to see that the filtered Boolean 
power of $F_i$ with respect to $B_i /I$ and 
the corresponding "quotient" of $\tau_i$ 
is a product of finite direct powers of subfields 
from $Rng(\tau_i )$. 
As in the proof of Lemma \ref{prod} we see 
that this product is a homomorphic image of 
$C(X_i ,F_i ,\tau_i )$. 
By the ring theory version of Proposition 5.1 from 
\cite{apps2} (see Introduction) any homomoprphic 
image of $C(X_i ,F_i ,\tau_i )$ is of this form. 
Thus the inverse limit of all finite quotients of 
$C(X_i ,F_i ,\tau_i )$ is the direct 
product of $F^{\omega}_i$ and finitely many direct 
powers $\tau_i (C)^{\omega}$, $C\in Rng(\tau_i )$.    

From this description we see that $(V,\oplus ,\otimes )$ 
is the direct product of finitely many direct powers 
$R^{\omega}$, where $R$ is a subfield of some $F_i$ as above. 
This is a contradiction. 
To see this precisely note that $R^{\omega}$ 
cannot be $\omega$-categorical. 
Let $f_1 , f_2 , ..., f_{i}, ... $ be a sequence 
of elements of $R^{\omega}$ with finite support 
such that for any pair of distinct $i,j$ we have
$|supp(f_i )|\not= |supp(f_j )|$. 
Then all ideals $f_i R^{\omega}$ are of pairwise distinct size. 
We see that for all $i\not= j$, the elements $f_i$ and $f_j$ 
have distinct types, which contradicts $\omega$-categoricity. 
$\Box$

\bigskip 

Let $G$ be an $\omega$-categorical group having a normal 
abelian subgroup $A$ of finite index. 
Consider $A$ as a $G/A$-module. 
Then $A$ is decomposed into a direct sum of finitely 
many modules $A_i$ such that each $A_i$ is a direct 
sum of $\omega$ copies of a finite indecomposable 
$G/A$-module (see \cite{baur} and Sections 11,12 and 
Apendix of \cite{BCM}). 
It is clear that the profinite completion of 
the module $A$ is an elementary extension of $A$.  
This implies that 
\begin{quote} 
when $G$ is an $\omega$-categorical abelian-by-finite 
group, $G$ is residually finite and is an elementary 
substructure of the completion $\overline{G}$. 
\end{quote} 
In particular under the BCM-conjecture 
\footnote{an $\omega$-categorical stable group is 
abelian-by-finite \cite{BCM}} 
we have that a stable $\omega$-categorical group 
is residually finite and is an elementary 
substructure of its profinite completion.    
Krzysztof Krupi\'{n}ski has suggested that 
this statement can be proved without the BCM-conjecture 
(and then this can be considered as a small 
confirmation of the BCM-conjecture). 
We think that an appropriate development of the 
methods of \cite{m1} allows a stronger statement: 
an $\omega$-categorical group without 
the strict order property is residually finite 
and is an elementary substructure of its profinite 
completion. 
Since the technique of \cite{m1} is quite involved 
we postpone this for a separate reseach.

\section{Measuring the set of commuting pairs}

Let $G$ be an $\omega$-categorical group, $H$ be 
an inert residually finite subgroup and $\overline{G}$ 
be the completion of $G$ with respect to 
$X= \{ gK : K$ is commensurable with $H$ and $g\in G\}$.   
Let $\mu$ be the Haar measure on $\overline{G}$. 
In the case when $\overline{G}$ is profinite and $\mu$ 
is normalized it follows from the main result of \cite{LP} 
that if there is a real number $\varepsilon >0$ such that 
$(\mu\times\mu ) (\{ (x,y)\in \overline{G}\times \overline{G} :xy=yx\})=\varepsilon$, 
then $\overline{G}$ is abelian-by-finite. 
In this section we consider some related conditions 
which work in the case of $\overline{G}$ without 
the assumption that $\overline{G}$ is profinite. 
We also describe some approximating functions which 
naturally arise in these respects. \parskip0pt  

Let $Y\subset X$ consist of finitely many $G$-orbits. 
By $G_{Y}$ ($\overline{G}_{Y}$ resp.) we denote 
the pointwise stabilizer of $Y$ in $G$ (in $\overline{G}$).   
We denote $G(Y):=G/G_Y$ and $\overline{G}(Y):=\overline{G}/\overline{G}_Y$. 
These groups acts faithfully on $Y$ and 
$\overline{G}(Y)$ is the completion of $G(Y)$ under 
the topology of this action. 
As in Proposition \ref{Aext} we see that every 
finitely generated subgroup of $\overline{G}(Y)$ is 
a homomorphic image of a finite subgroup of $G(Y)$. \parskip0pt 
 
Let $\rho_{com}(i)$ be the function which associates 
to each natural number $i$ the maximal number $n$ such 
that if $L$ is a subgroup of some $G(Y)$ where $Y\subset X$  
consists of finitely many $G$-orbits, and $|L|=i$, then 
$n \le |\{ (x,y)\in L\times L: xy=yx\}|$. 
If no $G(Y)$ has subgroups of size $i$, then 
we put $\rho_{com}(i)=\infty$. \parskip0pt 

In the case when $G$ is residually finite and $X$ consists 
of cosets of subgroups of finite index, the function 
$\rho_{com}(i)$ is closely connected with the value 
$(\mu\times\mu )(\{ (x,y)\in \overline{G}\times \overline{G} :xy=yx\})$. 
On the other hand it is clear that the square 
$\mu\times\mu$ is equal to the normalized Haar 
measure on $\overline{G}\times\overline{G}$. 
We will denote it by the same letter $\mu$. 
By a straightforward generalization of 
(the first proof of) Theorem 1(iii) from 
\cite{LP} we have the following statement. 

\begin{lem} 
Let $C$ be a compact group. 
If for some positive constant $\varepsilon \in \mathbf{R}$,  
$\mu_{C\times C} (\{ (x,y)\in C\times C :xy=yx\})>\varepsilon$, 
then $C$ is finite-by-abelian-by-finite. 
\footnote{in \cite{LP} this statement is formulated for 
profinite groups but with a slightly stronger conclusion: 
$C$ is abelian-by-finite}  
\end{lem} 

{\em Proof.} 
The proof of Lemma 5 from \cite{LP} without the 
last two sentences, works for any compact group 
\footnote{in \cite{LP} it concerns profinite groups}  
and shows that if the FC-centre of $C$ (i.e. the set 
of elements of $C$ having finitely many conjugates) is 
of finite index in $C$ then $C$ is finite-by-abelian-by-finite. 
Now the first proof of Theorem 1(iii) \cite{LP} works 
in our case. 
$\Box$

\bigskip 

Below we also use another observation from \cite{LP} 
(see Theorem 1(ii) and its proofs).  
This is the fact that in the profinite case the condition 
$\mu (\{ (x,y)\in C\times C :xy=yx\})=\varepsilon$ 
is equivalent to the condition that there are natural 
numbers $n_1$ and $n_2$ such that for any finite quotient 
$L$ of $C$ there is a subgroup $N<L$ with $|L:N|\le n_1$ 
and $|[N,N]|\le n_2$. 

\begin{lem} \label{comlem} 
Let $G$ be a residually finite group. 
Assume that there is a positive real number 
$\varepsilon <1$ such that for any finite quotient $L$ 
of $G$ we have $\rho_{com}(|L|) \ge \varepsilon |L|^2$. 
Then both $G$ and the profinite completion 
$\overline{G}$ are abelian-by-finite and there 
is an open normal subgroup $N < \overline{G}$ of 
nilpotency class two such that $|[N,N]|$ and 
$|\overline{G}/N|$ are both $\varepsilon$ bounded. 
\parskip0pt 

In particular if $G$ is finite-by-abelian-by-finite 
(for example, if $Th(G)$ is $\omega$-categorical and 
supersimple \cite{EW}), then both $G$ and $\overline{G}$ 
are abelian-by-finite. 
\end{lem}

{\em Proof.} 
Assume that for any finite quotient $L$ of $G$, 
$|\{ (x,y)\in L\times L: xy=yx\}|\ge \varepsilon |L|^2$. 
Since the set of commuting pairs is closed in 
$\overline{G}\times\overline{G}$ we easily see that 
the $\mu$-measure of this set is at least $\varepsilon$. 
By Theorem 1 of \cite{LP} we obtain that $\overline{G}$ is 
abelian-by-finite, and there is an open normal subgroup 
$N < \overline{G}$ of nilpotency class two such that $|[N,N]|$ 
and $|\overline{G}/N|$ are both $\varepsilon$ bounded. 
Since $G$ embeds into $\overline{G}$, $G$ is abelian-by-finite too. 

To see the second statement of the lemma note that the 
assumptions of this statement imply that there is a real number 
$\varepsilon <1$ such that for any finite quotient $L$ of 
$G$ there are normal subgroups $K<N<L$ such that $N/K$ is 
abelian and $|K|\cdot |L:N|^2 \le \varepsilon^{-1}$.   
Applying the inequality from \cite{pn} (see Section 1) 
$$
|\{ (x,y)\in L\times L: xy=yx\}|\ge |L|^2 /(|K|\cdot |L:N|^2 ) 
$$
we see that 
$|\{ (x,y)\in L\times L: xy=yx\}|\ge \varepsilon |L|^2$. 
Thus we can apply the basic statement of the lemma. 
$\Box$ 

\bigskip 

We now consider the case when $G$ is locally finite 
and $\overline{G}$ is locally compact. 
The following proposition is the main result 
of the section. 

\begin{prop} \label{comm} 
Let $G$ be a locally finite group having finite exponent 
and infinite inert residually finite subgroups 
(for example $G$ is $\omega$-categorical). 
Let $X$ be defined as above for some infinite 
inert residually finite $H\le G$. 
The following conditions are equivalent. \\  
(1) There is a positive real number $\varepsilon <1$ 
such that for any finite subgroup 
$L$ of $G(Y)$ with $Y\subset X$ consisting of finitely 
many $G$-orbits, we have 
$\rho_{com}(|L|) \ge \varepsilon |L|^2$; \\ 
(2) There is a positive real number $\varepsilon <1$  
such that for any infinite compact subgroup 
$C < \overline{G}$ we have   
$\mu_C (\{ (x,y)\in C\times C :xy=yx\})\ge\varepsilon$, 
where $\mu_C$ is the normalized Haar measure on $C$; \\
(3) Both $G$ and $\overline{G}$ are finite-by-abelian-by-finite. 

Each of the conditions (1), (2) and (3) implies that 
every infinite inert residually finite subgroup $K<G$  
(and the completion $\overline{K}$) is abelian-by-finite. 

The group $G$ is abelian-by-finite if and only if 
there is a number $l$ such that every infinite inert 
residually finite subgroup $K<G$  
(and the completion $\overline{K}$) is abelian-by-$\le l$. 
\end{prop}

{\em Proof.} 
(1) $\Rightarrow$ (3). 
We will use the following theorem of P.Neumann. 
For each positive number $\varepsilon$ there are natural numbers 
$n_1$ and $n_2$ such that if a finite group $F$ has at least 
$\varepsilon |F|^2$ commuting pairs, then there is a normal 
subgroup $N<F$ with $|F:N|\le n_1$ and $|[N,N]|\le n_2$ \cite{pn}.   
Moreover in this case $n^2_1 n_2 \le \varepsilon^{-1}$. 
\parskip0pt 

Let $Y\subset X$ consist of finitely many $G$-orbits 
and $L<G(Y)$ be finite. 
Then by Neumann's theorem there is a normal subgroup 
$N < L$ such that $|L:N|\le n_1$ and $|[N,N]|\le n_2$ 
where $n_1 ,n_2$ are determined by $\varepsilon$ as above. 
\parskip0pt 

Now note that this statement holds for any 
finite subgroup of $G$. 
This follows from the fact that if $F$ is a finite 
subgroup of $G$, then $F$ can be realized as 
a subgroup of $G(Y)$ for an appropriate $Y$ as above.  
Indeed since $H$ is inert, $F$ normalizes a subgroup 
$H_1 <H$ of finite index in $H$ such that 
$F\cap H_1 =\{ 1\}$. 
Then for $Y=\{ gH_1 :g\in G\}$ (which is a $G$-orbit) 
there is an embedding of $F$ into $G(Y)$. 
Thus we can find a normal subgroup $N < F$ such 
that $|F:N|\le n_1$ and $|[N,N]|\le n_2$. \parskip0pt 

We present $G$ as the union of an $\omega$-chain 
of finite groups $L_0 <L_1 <...<L_i <...$. 
For each $i$ we find a normal $N_i <L_i$ with 
$|L_i :N_i |\le n_1$ and $|[N_i ,N_i ]|\le n_2$.  
Choosing a subsequence if necessary we can arrange 
that for any pair $g_1 ,g_2 \in G$ 
there is a natural number $j$ such that we have 
$(\forall i\ge j) L_i \models (g_1 N_i =g_2 N_i )$ or 
$(\forall i\ge j) L_i \models (g_1 N_i \not=g_2 N_i )$. 
As a result we see that $N_{\omega}:=\{g\in G: g\in N_i$ 
for infinitely many $i \}$ is a subgroup of index $\le n_1$ 
in $G$. 
Since each $|[N_i ,N_i ]|\le n_2$ we see that 
$[N_{\omega},N_{\omega}]\le n_2$. 
This obviously means that $\overline{N}_{\omega}$ 
is finite-by-abelian. 
Since $\overline{N}_{\omega}$ is of index $\le n_1$ 
in $\overline{G}$, this proves (3). 
\parskip0pt 

Let us prove (3) $\Rightarrow$ (2).  
Find $n_1$ and $n_2$ such that there is a normal subgroup 
$N<G$ of index $\le n_1$ such that $|[N,N]|\le n_2$. 
Let $C$ be a compact subgroup of $\overline{G}$. 
Then $C$ and all its finite quotients are 
$(\le n_2)$-by-abelian-by-$(\le n_1 )$. 
This implies that there is a positive real 
number $\varepsilon <1$ such that for any 
finite quotient $L$ of $C$ we have  
$|\{ (x,y)\in L\times L: xy=yx\}|\ge \varepsilon |L|^2$  
(by an easy inequality from Section 1 of \cite{pn} 
we may put $\varepsilon = (n^2_1 n_2)^{-1}$). 
This implies (2). 

(2) $\rightarrow$ (1). 
Note that if $H_0$ is a residually finite inert 
subgroup of $G$ which is commensurable with $H$, 
then the closure of $H_0$ in $\overline{G}$ is profinite. 
Condition (2) implies that $\overline{H}_0$ 
satisfies the assumption of Theorem 1 of \cite{LP}. 
Thus $\overline{H}_0$ has an open normal subgroup 
$N < \overline{H}_0$ of nilpotency class two such 
that $|[N ,N ]|$ and $|\overline{H}_0 /N |$ 
are both $\varepsilon$ bounded. 
\parskip0pt 

On the other hand if $F$ is a finite subgroup of $G$, 
then $F$ can be realized as a subgroup of 
an inert residually finite $H_0$ which is 
commensurable with $H$.  
Indeed since $H$ is inert, $F$ normalizes 
a subgroup $H_1 <H$ of finite index in $H$. 
Then $\langle H_1 ,F\rangle$ works as $H_0$. 
The fact that $\langle H_1 ,F\rangle$ is inert 
and commensurable with $H$ is obvious. 
To see residual finiteness note that if $H_2$ 
is a normal subgroup of $H_1$ of finite index, 
then it is of finite index in $\langle H_1 ,F\rangle$ 
and thus contains a normal subgroup of 
$\langle H_1 ,F\rangle$ of finite index.  
Since $H_1$ is residually finite, 
$\langle H_1 ,F\rangle$ is residually finite too. 
\parskip0pt 

Thus as in the proof of Lemma \ref{comlem} for 
every finite subgroup $F<G$ we can find 
a normal subgroup $N < F$ such that both $|[N,N]|$ 
and $|F:N|$ are $\varepsilon$ bounded. 
By an argument from \cite{pn} (or the proof of 
Lemma \ref{comlem}) for any finite subgroup $L$ of 
$G(Y)$ with $Y\subset X$ consisting of finitely many $G$-orbits, 
we have (correcting $\varepsilon$ if necessary) 
$\rho_{com}(|L|) \ge \varepsilon |L|^2$.  
\parskip0pt 

The statement of the proposition that all inert 
residually finite subgroups are abelian-by-finite, 
follows from Lemma \ref{comlem}. 

To see the last statement note that the condition that 
there is a number $l$ such that every infinite inert 
residually finite subgroup $K<G$  
(and the completion $\bar{K}$) is abelian-by-$\le l$
implies that every finite subgroup of $G$ is 
abelian-by-$\le l$ (see the proof of (2) $\Rightarrow$ (1)). 
Now applying the argument of (1) $\Rightarrow$ (3) 
we see that $G$ is abelian-by-$\le l$. 
$\Box$ 

\bigskip 

The function $\rho_{com}$ is connected with 
some other functions naturally arising in these respects. 
For example consider the following asymptotic condition. 
Let $\rho_{r}(i)$ be the function (defined by $G$) which 
associates to each natural number $i$ the minimal number 
$m$ such that whenever $L$ is a quotient of $G$ with 
$|L|=i$, then $H/H_{L}$ has rank at most $m$ for 
all subgroups $H$ of $L$. 
Here $H_L$ is the maximal normal subgroup of $L$ 
contained in $H$ ({\em the core of $H$ in $L$}).  
It is worth mentioning that Proposition 4.3 of 
\cite{KhukhroSmith} states that there is a function 
$\beta :\omega\rightarrow\omega$ such that 
if $L$ as above is a finite group with $H/H_L$ of 
rank at most $r$ for all $H<L$, then $L$ is 
abelian-by-(group of rank $<\beta (r)$). 
Assuming that $G$ is locally finite of exponent 
$m\in\omega$ we can apply Zelmanov's solution of 
the restricted Burnside problem to derive that in fact 
the function $\beta$ can be chosen so that 
$L$ is abelian-by-(group of order $<\beta (r)$).  
\parskip0pt 

Another interesting function can be defined  
using the approach of Baudisch (in \cite{Baudisch}) 
to the BCM-conjecture. 
Assume that $G$ is additionally a two-step-nilpotent 
group and $U:=G/[G,G]$. 
In the situation of \cite{Baudisch} it is also assumed 
that $U$ and $W:=[G,G]$ are vector spaces over $GF(p)$ 
(this case is basic for the BCM-conjecture).     
Then the commutator $[x,y]$ becomes an alternating 
bilinear map $U\times U \rightarrow W$ over $GF(p)$.  
If $L$ is a quotient of $G$ of size $i$, then there is 
a vector space homomorphism $f_L$ from the exterior 
square $\Lambda^{2}(L/(L\cap W))$ to $W\cap L$ 
\footnote{if $L=G/N$, then by $L\cap W$ we denote the elements of the form $wN$, $w\in W$}   
which makes commutative the diagram consisting of the maps 
$[,]: L/(L\cap W)\times L/(L\cap W) \rightarrow L\cap W$ and 
$\wedge : L/(L\cap W)\times L/(L\cap W) \rightarrow \Lambda^{2}(L/(L\cap W))$. 
By $\rho_{\wedge}(i)$ we denote the maximal $k$ such 
that $k\cdot dim(L/(L\cap W))\le dim(Ker(f_L ))$ for 
all $L$ with $|L|=i$.  
Baudisch studied in \cite{Baudisch} the related 
condition that there is a natural number $k$ such that 
$dim(Ker(f_H ))\le k\cdot dim(H/(H\cap W))$ for all 
finite subgroups $H$ of $G$ (in this case Baudisch 
says that $G$ {\em has few relations}).  
\parskip0pt 

The following computations demonstrate some connections 
of these functions with $\rho_{com}(i)$. 

\begin{lem} 
(1) If the function $\beta$ is chosen as above, then 
$\beta(\rho_{r}(i))^{-2}\cdot i^2 \le \rho_{com}(i)$. \\     
(2) In the situation of the function $\rho_{\wedge}$ 
we have 
$i^{(2\rho_{\wedge}(i)+1-log_p i)/2}\cdot i^2 \le \rho_{com}(i)$.  
\end{lem} 

{\em Proof.} 
By an argument from \cite{pn} (p.456) for any finite 
group $L$ with normal subgroups $K<N<L$ such that 
$N/K$ is abelian we have the following inequality:  
$|\{ (x,y)\in L\times L: xy=yx\}|\ge (|K|\cdot |L:N|^2)^{-1} |L|^2$. 
Assuming that the function $\beta$ is as above, 
this inequality shows that: 
$\beta(\rho_{r}(i))^{-2}\cdot i^2 \le \rho_{com}(i)$.    
\parskip0pt 

In the case of the function $\rho_{\wedge}$ we 
preserve the notation above. 
Then the argument of the previous paragraph shows 
that $\rho_{com}(i) \ge i^2 \cdot |W\cap L|^{-1}$. 
On the other hand when $L$ is a quotient of $G$ of size $i$ 
we have $|W\cap L|=p^{dim(W\cap L)}$, 
$|L/(W\cap L)|=p^{dim(L/(W\cap L))}$ and 
$dim(W\cap L) = dim(\Lambda^{2} L/(L\cap W)) - dim(Ker(f_L ))$. 
Since 
$dim(\Lambda^2 L/(L\cap W))=dim(L/(L\cap W))\cdot (dim(L/(L\cap W))-1)/2$ 
and 
$dim(Ker(f_L ))\ge \rho_{\wedge}(i)dim(L/(L\cap W))$,   
we have 
$dim(W\cap L)\le dim(L/(L\cap W))(dim(L/(L\cap W))-1-2\rho_{\wedge}(i))/2$ 
$\le (log_p i)(log_{p}i -1-2\rho_{\wedge}(i)) /2$.   
Thus $|W\cap L|^{-1} \ge p^{(log_p i)(2\rho_{\wedge}(i)+1-log_p i )/2}$, 
which implies the required inequality. 
$\Box$


\begin{thebibliography}{9} 
\bibitem{apps2} A.B.Apps. Boolean powers of groups. 
{\em Math. Proc. Cambridge Philos. Soc.}, 91(1982), 375 - 395. 
\bibitem{apps} A.B.Apps. On the structure of 
$\aleph_0$-categorical groups, {\em J. Algebra}, 81(1883), 320 - 339. 
\bibitem{archmac} R.Archer and D.Macpherson. Soluble 
$\omega$-categorical groups. {\em Math. Proc. Cambridge Philos. Soc.}, 
121(1997), 219 - 227. 
\bibitem{Baudisch} A.Baudisch. Closures in $\aleph_0$-categorical 
bilinear maps. {\em J.Symb.Logic}, 65(2000), 914 - 922. 
\bibitem{belyaev} V.V.Belyaev. Locally finite groups with finite 
non-separable subgroups. {\em Syberian Math.J.}, 34(1993), 23 - 41. 
\bibitem{baur} W.Baur. $\aleph_0$-categorical modules.
{\em J. Symb. Logic}, 40(1975), 213 - 226.
\bibitem{BCM} W.Baur, G.Cherlin and A.Macintyre. Totally categorical 
groups and rings. {\em J.Algebra}, 57(1979), 407 - 440. 
\bibitem{EW} D.M.Evans and F.Wagner. Supersimple $\omega$-categorical 
groups and theories. {\em J.Symb.Logic}, 65(2000), 767 - 775.  
\bibitem{gorenstein} Gorenstein. {\em Finite Groups} 
(Chelsea Publishing Company, New York, 1980). 
\bibitem{KhukhroSmith} E.I.Khukhro and H.Smith. Locally finite 
groups with all subgroups normal-by-(finite rank). {\em J.Algebra}, 
200(1998), 701 - 717. 
\bibitem{KoTits} S.Koppelberg and J.Tits. Une propriete des produits 
directs infinis de groupes finis isomorphes. 
{\em C. R. Math. Acad. Sci. Paris}, Ser. A 279 (1974), 583 - 585. 
\bibitem{LP} L.Levai and L.Pyber. Profinite groups with many 
commuting pairs of involutions. {\em Arch. Math.}, 75(2000), 1 - 7. 
\bibitem{macros} A.J.Macinryre and J.G.Rosenstein. 
$\aleph_0$-categoricity for rings without nilpotent elements and 
for boolean structures. {\em J.Algebra}, 43(1998), 483 - 500. 
\bibitem{m1} H.D. Macpherson. Absolutely ubiquitous
structures and $\omega$-categorical groups. {\em Quart. J. Math.}
Oxford (2), 39(1988), 483-500.
\bibitem{neumann} P.Neumann. On the structure of standard wreath 
products of groups. {\em Math. Z.} 84(1964), 343 - 373. 
\bibitem{pn} P.Neumann. Two combinatorial problems in group theory. 
{\em Bull.London Math.Soc.}, 21(1989), 456 - 458.  
\bibitem{newelski} L.Newelski. Small profinite groups. {\em J.Symb.Logic}, 
66(2001), 859 - 872.
\bibitem{rosendal} Ch. Rosendal. A topological version of 
the Bergman property. ArXiv:math.LO/0509670 v1 (28 Sep 2005), 33 pps. 
\bibitem{ros} J.G.Rosenstein. $\aleph_0$-Categoricity of Groups. 
{\em J.Algebra}, 25(1973), 435 - 467. 
\bibitem{sst} J.Saxl, S.Shelah and S.Thomas. Infinite products 
of finite simple groups. {\em Trans. Amer. Math. Soc.}, 
348(1996), 4611 - 4641. 
\bibitem{wilson} J.Wilson. The algebraic structure of 
$\aleph_0$-categorical groups. In: {\em Groups St.Andrews 1981}, 
London Math. Soc. Lecture Notes no. 71, Eds. C.M.Campbell, 
E.F.Robertson (Cambridge University Press, 1983), pps. 345 - 358.  
\bibitem{WilsonProf} J.Wilson. {\em Profinite Groups}. 
(Oxford University Press, Oxford, 1998).  
\end{thebibliography}
\end{document}